\documentclass[11pt,oneside]{amsart}

\usepackage{bm}

\usepackage{fix-cm}
\DeclareMathAlphabet{\mathcal}{OMS}{cmsy}{m}{n}

\usepackage{tikz-cd}

\usepackage{tcolorbox}

\usepackage[arrow,curve,matrix,tips,2cell]{xy}

\usepackage{fancyhdr}

\usepackage{anyfontsize}
\usepackage{mathrsfs}

\usepackage{accents}

\usepackage{bbm}

\usepackage{amsmath}
\usepackage{amscd}
\usepackage{amssymb}
\usepackage{latexsym}
\usepackage{url}

\usepackage{enumitem}

\usepackage{graphicx}

\newtheorem{theorem}{Theorem}[section]
\newtheorem*{theorem*}{Theorem}
\newtheorem{lemma}[theorem]{Lemma}
\newtheorem*{lemma*}{Lemma}
\newtheorem{corollary}[theorem]{Corollary}
\newtheorem{proposition}[theorem]{Proposition}

\newtheorem{remark}[theorem]{Remark}
\newtheorem{definition}[theorem]{Definition}
\newtheorem*{definition*}{Definition}

\newtheorem{question}[theorem]{Question}
\newtheorem*{question*}{Question}

\newtheorem{example}[theorem]{Example}
\newtheorem{examples}[theorem]{Examples}

\makeatletter
\def\revddots{\mathinner{\mkern1mu\raise\p@
\vbox{\kern7\p@\hbox{.}}\mkern2mu
\raise4\p@\hbox{.}\mkern2mu\raise7\p@\hbox{.}\mkern1mu}}
\makeatother 
\newcommand{\bgl}{\begin{equation}} 
\newcommand{\egl}{\end{equation}}
\newcommand{\bgloz}{\begin{equation*}} 
\newcommand{\egloz}{\end{equation*}}
\newcommand{\bgln}{\begin{eqnarray}} 
\newcommand{\egln}{\end{eqnarray}}
\newcommand{\bglnoz}{\begin{eqnarray*}} 
\newcommand{\eglnoz}{\end{eqnarray*}}
\newcommand{\btheo}{\begin{theorem}}
\newcommand{\etheo}{\end{theorem}}
\newcommand{\btheooz}{\begin{theorem*}}
\newcommand{\etheooz}{\end{theorem*}}

\newcommand{\blemma}{\begin{lemma}}
\newcommand{\elemma}{\end{lemma}}
\newcommand{\blemmaoz}{\begin{lemma*}}
\newcommand{\elemmaoz}{\end{lemma*}}
\newcommand{\bproof}{\begin{proof}}
\newcommand{\eproof}{\end{proof}}
\newcommand{\bbew}{\begin{beweis}}
\newcommand{\ebew}{\end{beweis}}
\newcommand{\bremark}{\begin{remark}\em}
\newcommand{\eremark}{\end{remark}}
\newcommand{\bdefin}{\begin{definition}}
\newcommand{\edefin}{\end{definition}}
\newcommand{\bdefinoz}{\begin{definition*}}
\newcommand{\edefinoz}{\end{definition*}}
\newcommand{\bex}{\begin{example}}
\newcommand{\eex}{\end{example}}
\newcommand{\bexs}{\begin{examples}}
\newcommand{\eexs}{\end{examples}}

\newcommand{\bprop}{\begin{proposition}}
\newcommand{\eprop}{\end{proposition}}
\newcommand{\bcor}{\begin{corollary}}
\newcommand{\ecor}{\end{corollary}}
\newcommand{\bfa}{\begin{cases}} 
\newcommand{\efa}{\end{cases}}

\newcommand{\bquestion}{\begin{question}}
\newcommand{\equestion}{\end{question}}
\newcommand{\bquestionoz}{\begin{question*}}
\newcommand{\equestionoz}{\end{question*}}
\newcommand{\cO}{\mathcal O}

\newcommand{\cZ}{\mathcal Z}

\def\1z{\mathbb{1}}
\newcommand{\onto}{\twoheadrightarrow} 

\def\SEMI{\mbox{$\times\kern-2pt\vrule height5pt width.6pt \kern3pt $}}

\newcommand{\norm}[1]{\left\|#1\right\|} 
\newcommand{\defeq}{\mathrel{:=}} 

\newcommand{\dop}{\text{: }} 
\newcommand{\lge}{\left\{} 
\newcommand{\rge}{\right\}} 
\newcommand{\gekl}[1]{\lge #1 \rge} 
\newcommand{\menge}[2]{\gekl{ #1 \dop #2 }} 
\newcommand{\oset}[2]{%
  \mathop{#2}\limits^{\vbox to -1.66ex{%
  \kern -1.4ex\hbox{$#1$}\vss}}}

\newcommand{\pari}{\setlength{\parindent}{0.5cm} \setlength{\parskip}{0cm}}
\newcommand{\nopar}{\setlength{\parindent}{0cm} \setlength{\parskip}{0cm}}

\usepackage{newtxtext}
\usepackage{newtxmath}

\begin{document}

\title{C$^*$-diagonals in AH-algebras arising from generalized diagonal connecting maps: spectrum and uniqueness}

\thispagestyle{fancy}

\author{Ali I. Raad}

\address{Ali I. Raad, Mathematics and Science Department, American University in Bulgaria, Office 307 BAC, ul. Svoboda Bacharova 8, 2700 Blagoevgrad, Bulgaria}
\email{araad@aubg.edu}


\thanks{This project has been supported by the Internal KU Leuven BOF project C14/19/088 and project G085020N funded by the Research Foundation Flanders (FWO)}

\begin{abstract}
We associate a Bratteli-type diagram to AH-algebras arising from generalized diagonal connecting maps. We use this diagram to give an explicit description of the connected components of the spectrum of an associated canonical C$^*$-diagonal. We introduce a topological notion on these connected components, that of being spectrally incomplete, and use it as a tool to show how various classes of AI-algebras, including certain Goodearl algebras and AH-algebra models for dynamical systems $([0,1],\sigma)$, do not admit unique inductive limit Cartan subalgebras. We focus on a class of spectrally complete C$^*$-algebras, namely the AF-algebras, and discuss the uniqueness of their inductive limit Cartan subalgebras.
\end{abstract}

\maketitle


\setlength{\parindent}{0cm} \setlength{\parskip}{0.5cm}

\section{Introduction}
In recent years there has been a growing interest in the study of Cartan subalgebras of C$^*$-algebras. In \cite{Li16,Li_DQH} Li obtained strong rigidity results for crossed product C$^*$-algebras arising from dynamical systems as well as for coarse geometric structures arising from group dynamics, all via Cartan subalgebras. Thus it was shown that Cartan subalgebras establish a link between C$^*$-algebras on the one hand, and dynamical systems and geometric group theory on the other. 

Li has also established a link between Cartan subalgebras and the classification programme for C$^*$-algebras. This programme, due to many hands (see for example \cite{KP, Phi, GLNa, GLNb, EGLN, TWW, EGLN17a,EGLN17b,GLI,GLIII}) aims to classify C$^*$-algebras by $K$-theoretic and tracial data, and has seen tremendous progress over the years. In \cite{Li18} it is proved that every classifiable C$^*$-algebra has a Cartan subalgebra. 

Cartan subalgebras are also intimately linked to topological groupoids. Many important classes of C$^*$-algebras can be described naturally as groupoid C$^*$-algebras, for instance C$^*$-algebras attached to group actions (crossed products) or semigroup actions on topological spaces (Cuntz-Krieger algebras or graph algebras) or AF-algebras. By the work of Kumjian \cite{Kum} and Renault \cite{Ren} it was shown that separable C$^*$-algebras have Cartan subalgebras if and only if they arise as twisted groupoid C$^*$-algebras. These twisted groupoids are of a specific nature: namely they are locally compact, second-countable, topologically principal, Hausdorff and étale. The result was then extended to all C$^*$-algebras (i.e. even non-separable ones) by Kwa\'sniewski and Meyer \cite{KM}, and independently by the author \cite{Raa}. Here the conditions on the twisted groupoid weaken to being effective, locally compact, Hausdorff and étale.

All these motivating factors led to interest in existence and uniqueness questions for Cartan subalgebras. Li and Renault initiated the systematic study of existence and uniqueness of Cartan subalgebras in C$^*$-algebras in \cite{LR19}. Barlak and Raum classified Cartan subalgebras of dimension drop algebras with coprime parameters in \cite{bar-raum}. White and Willett studied the question of uniqueness of Cartan subalgebras in uniform Roe algebras in \cite{white}. Li and the author constructed C$^*$-diagonals in \cite{LR22}, which extended the work of \cite{Li18} by finding Cartan subalgebras in many classes of non-classifiable C$^*$-algebras. 

In this work we aim to initiate the question of uniqueness of Cartan subalgebras in AH-algebras arising from generalized diagonal connecting maps, which was not covered in \cite{LR22}. Section 2 will provide some preliminaries on topological groupoids and inductive limit Cartan subalgebras, AH-algebras which arise from generalized diagonal connecting maps, as well as how the associated C$^*$-diagonal is obtained, as is done in \cite{LR22}. In Section 3 we associate a specific Bratteli diagram $\mathcal{B}_A$ to AH-algebras $A$ arising from generalized diagonal connecting maps, by labeling edges via the eigenvalue functions that correspond to the connecting map. We then show how the spectrum of the associated C$^*$-diagonal $C_A$ can be fully described through $\mathcal{B}_A$. More specifically, we obtain:

\textbf{Theorem A} (cf. Theorem \ref{theorem:connected components}). \emph{Let $(A,C_A)$ be a pair of AH-algebra $A$ arising from generalized diagonal connecting maps and associated C$^*$-diagonal $C_A$, with associated Bratteli diagram $\mathcal{B}_A$. Then the connected components of the spectrum of $C_A$ are precisely the inverse limits over infinite paths of $\mathcal{B}_A$.}

We then define a topological notion on the connected components of the spectrum of the associated C$^*$-diagonal $C_A$. This notion, termed \emph{spectrally incomplete}, captures the notion that the connected components of the spectrum do not exhaust (up to homeomorphism) all possible inverse limits over the base spaces that give rise to the AH-algebra. We relate this notion to uniqueness of inductive limit Cartan subalgebras. In Section 4 we prove that if an AI-algebra is spectrally incomplete then it does not admit unique inductive limit Cartan subalgebras:

\textbf{Theorem B} (cf. Theorem \ref{theorem:AI non uniqueness}). 
\emph{Simple and unital spectrally incomplete AI-algebras arising from interval algebra building blocks do not have unique inductive limit Cartan subalgebras.}

We use this theorem to then show that certain Goodearl algebras and AH-algebra models for dynamical systems do not have unique inductive limit Cartan subalgebras:

\textbf{Theorem C} (cf. Corollaries \ref{cor: Goodearl} and \ref{cor: AH algebra models}). \emph{Goodearl algebras defined over unit interval base spaces do not admit unique inductive limit Cartan subalgebras, neither do AH-algebra models for a dynamical system $([0,1],\sigma)$.}

In Section 5 we consider the class of AF-algebras, as these are spectrally complete. We adapt the ideas of Elliott's original proof for classification of AF-algebras via the $K_0$-group \cite{ellaf} to show that the dimension group is a complete invariant for inductive limit Cartan subalgebras:

\textbf{Theorem D} (cf. Theorem \ref{main theorem}) \emph{Let $A$ and $B$ be unital AF-algebras with inductive limit Cartan subalgebras $C$ and $D$, respectively. Assume there exists an isomorphism
$$\alpha: (K_{0}(A),K_{0}(A)^+,[1_A]_0) \rightarrow (K_{0}(B),K_{0}(B)^+,[1_B]_0).$$
Then there exists a $*$-isomorphism $\phi:A \rightarrow B$ such that $K_0(\phi)=\alpha$ and $\phi(C) = D$. In particular, AF-algebras have unique inductive limit Cartan subalgebras. }

Although this last result can be obtained by for instance combining \cite[Theorem~5.7]{Pow} together with Elliott's classification, it has come to our attention that not many mathematicians who are working within the field of Cartan subalgebras of C$^*$-algebras are aware of it, as the  references which obtain this result precede the current language used. We proved Theorem D before being made aware that it could be obtained by these means. Hence it would be important to re-emphasize such a result given the recent growing interest for Cartan subalgebras. In addition, our techniques of proof are different and are more in line with the techniques used to prove Theorem B.

\textbf{Acknowledgements.} Parts of this work is based on the author's PhD results in \cite{Raa_PhD}, and therefore the author would like to thank Xin Li for helpful comments during that period. The author would also like to thank the anonymous referee for their thorough reading and provision of useful comments.
\section{Preliminaries}

\subsection{Groupoids and Inductive Limit Cartan Subalgebras}

Here we provide a brief and concise introduction to the objects that are of interest for this paper. For further details on groupoids and their C$^*$-algebras, one may further consult \cite{Sim} amongst other references. For further details on twisted groupoids and Cartan subalgebras, one may further consult \cite{Ren} amongst other references. For further details on inductive limit Cartan subalgebras, please consult Section 2.3.3 in \cite{Raa}, Section 3 in \cite{BL17} and Section 5 in \cite{Li18}.

A \emph{groupoid} $\mathcal{G}$ is a non-empty set together with a set $\mathcal{G}^{(2)}\subseteq \mathcal{G}\times \mathcal{G}$ consisting of \emph{composable pairs} upon which a multiplication map $m:\mathcal{G}^{(2)}\rightarrow \mathcal{G}$ is defined. A groupoid also has an involutive map $^{-1}:\mathcal{G}\rightarrow\mathcal{G}$ which acts like an inverse. One assumes that $m$ is associative wherever it is defined, and that $(g,g^{-1})$ and $(g^{-1},g)$ are composable pairs for all $g\in \mathcal{G}$. In this way a groupoid generalizes the notion of a group as multiplication may not necessarily be defined everywhere.  

The \emph{unit space} of $\mathcal{G}$ is $\mathcal{G}^0=\{gg^{-1} : g\in\mathcal{G}\}$. Associated to a groupoid is a \emph{source map} and \emph{range map}, $s$ and $r$, respectively. They are maps $\mathcal{G}\rightarrow\mathcal{G}^0$ defined by $s(g)=g^{-1}g$ and $r(g)=gg^{-1}$. A \emph{topological groupoid} is a groupoid $\mathcal{G}$ endowed with a topology making $m$ and $^{-1}$ continuous ($\mathcal{G}^{(2)}\subseteq\mathcal{G}\times \mathcal{G}$ is endowed with the subspace topology of the product topology). A topological groupoid is called \emph{étale} if the source and range maps are local homeomorphisms (where $\mathcal{G}^0$ is endowed with the subspace topology). Almost all the groupoids we will consider in this paper are étale and Hausdorff. This will imply, amongst other things, that the unit space is clopen, and that \emph{source fibers} and \emph{range fibers}, $\mathcal{G}_x=\{g \in \mathcal{G} : s(g)=x\}$ and $\mathcal{G}^y=\{g\in \mathcal{G} : r(g)=y\}$, are discrete, respectively. The \emph{isotropy fiber at $x$} is $G^x_x=\{g\in \mathcal{G} : s(g)=r(g)=x\}$, and the \emph{isotropy bundle} is $\mathcal{G}^\prime = \bigcup\limits_{x\in \mathcal{G}^0}\mathcal{G}^x_x=\{g\in\mathcal{G} : s(g)=r(g)\}$. A groupoid is called \emph{principal} if $\mathcal{G}^\prime=\mathcal{G}^0$. It is \emph{topologically principal} if the elements of $\mathcal{G}^0$ that have trivial isotropy are dense in $\mathcal{G}^0$. It is \emph{effective} if $\mathrm{int}(\mathcal{G}^\prime)=\mathcal{G}^0.$ 

Almost all the groupoids in this paper will also be locally compact. This allows us to build C$^*$-algebras out of them. A \emph{twisted groupoid} $(\mathcal{G},\Sigma)$ consists of a Hausdorff étale groupoid $\mathcal{G}$ together with a Hausdorff, étale and locally compact groupoid $\Sigma$, called the \emph{twist}, such that $$\mathcal{G}^0\times \mathbb{T} \hookrightarrow \Sigma \twoheadrightarrow \mathcal{G}$$is a central groupoid extension. This implies that $\Sigma$ is a principal circle bundle and we may identify $\mathcal{G}$ with $\Sigma/\mathbb{T}$.

Whenever we say that a twisted groupoid $(\mathcal{G},\Sigma)$ has a certain property, it means that $\mathcal{G}$ has that property. From a twisted étale locally compact Hausdorff groupoid $(\mathcal{G},\Sigma)$ we can build a (reduced) \emph{twisted groupoid C$^*$-algebra} $C^*_r(\mathcal{G},\Sigma)$. We say that a C$^*$-algebra $A$ has a \emph{twisted groupoid model} if $A\cong C^*_r(\mathcal{G},\Sigma)$ for some twisted groupoid $(\mathcal{G},\Sigma)$.

We now state the definition of a Cartan subalgebra of a C$^*$-algebra. Note that we are omitting the condition of containing an approximate unit for the ambient C$^*$-algebra (which is present in the original definition), as this was shown to be redundant by Pitts (see \cite{Pit}):
\bdefin[\cite{Kum,Ren}]
\label{def:CartanSubalgebra}
A sub-C$^*$-algebra $C$ of a C$^*$-algebra $A$ is called a Cartan subalgebra if
\nopar

\begin{itemize}
\item $C$ is maximally Abelian;
\item $C$ is regular, i.e., $N_A(C) \defeq \menge{n \in A}{n C n^* \subseteq C \ \text{and} \ n^* C n \subseteq C}$ generates $A$ as a C$^*$-algebra;
\item there exists a faithful conditional expectation $P: \: A \onto C$.
\end{itemize}

A pair $(A,C)$, where $C$ is a Cartan subalgebra of a C$^*$-algebra $A$, is called a Cartan pair.
\pari

A Cartan subalgebra $C$ is called a C$^*$-diagonal if $(A,C)$ has the unique extension property, i.e., every pure state of $C$ extends uniquely to a (necessarily pure) state of $A$.
\edefin

By the work of Kumjian \cite{Kum} and Renault \cite{Ren}, every separable C$^*$-algebra with a Cartan subalgebras has a twisted groupoid model which is locally compact, second-countable, topologically principal, Hausdorff and étale. Kwa\'sniewski and Meyer \cite{KM}, and independently the author \cite{Raa}, extend this result to all C$^*$-algebras (i.e. the condition of separability is dropped). The twisted groupoid models in this context have weaker conditions; namely they are effective, locally compact, Hausdorff and étale. The other direction of this result also holds. Namely given an effective, locally compact, Hausdorff and étale twisted groupoid $(\mathcal{G},\Sigma)$, then $(C^*_r(\mathcal{G},\Sigma),C_0(\mathcal{G}^0))$ is a Cartan pair.
\par

Using this correspondence, Barlak and Li \cite{BL17} develop a machinery that builds Cartan subalgebras in inductive limit C$^*$-algebras. First they study this at the level of connecting maps. They show that a surjective twisted groupoid homomorphism $(H,T) \rightarrow (\mathcal{G},\Sigma)$ which admits certain technical conditions gives rise to a $*$-homomorphism $C^*_r(\mathcal{G},\Sigma)\rightarrow C^*_r(H,T)$, whilst an injective twisted groupoid homomorphism $(H,T) \rightarrow (\mathcal{G},\Sigma)$ which admits certain technical conditions gives rise to a $*$-homomorphism $C^*_r(H,T)\rightarrow C^*_r(\mathcal{G},\Sigma)$. A homomorphism between (twisted) groupoids that gives rise to a $*$-homomorphism between C$^*$-algebras is called a \emph{groupoid model} for the $*$-homomorphism. 

By combining this with the work of Li in \cite{Li18} it is proved that given an inductive limit C$^*$-algebra $A=\varinjlim(A_n,\phi_n)$ with Cartan subalgebras $C_n\subseteq A_n$, if the connecting maps $\phi_n$ fully preserve the structure of the Cartan subalgebras, then the groupoid models of the Cartan pairs, $(\mathcal{G}_n,\Sigma_n)$, admit certain morphisms that induce the connecting map morphisms. Furthermore, a limit twisted groupoid can be built out of the building block twisted groupoids and the groupoid models of the connecting maps, which possesses the correct characteristics to make it the twisted groupoid model for the Cartan pair $(\varinjlim(A_n,\phi_n),\varinjlim(C_n,\phi_n)).$

Finally, we will need the following definition in our paper:
\bdefin
Let $\mathfrak{c}$ be a class of C$^*$-algebras. For an inductive limit C$^*$-algebra $A=\varinjlim(A_n,\phi_n)$, where $A_n\in\mathfrak{c}$ , a Cartan subalgebra of the form $C=\varinjlim(C_n,\phi_n)$, with $\mathfrak{c} \ni C_n \subseteq A_n$, will be called a \emph{$\mathfrak{c}$-inductive limit Cartan subalgebra}, or an \emph{inductive limit Cartan subalgebra} if it is clear what the class $\mathfrak{c}$ is.
\edefin

\subsection{AH-algebras from Generalized Diagonal Connecting Maps}

We will consider AH-algebras $A=\varinjlim(A_n,\phi_n)$ of a certain form, where for each $n\in \mathbb{N}$ we have 
\begin{equation*}
    A_n=\bigoplus\limits_{i=1}^{N_n}A_n^i,\;\;     A_n^i=p_n^i(C(Z_n^i)\otimes M_{r_n^i})p_n^i,
\end{equation*}
where each $Z_n^i$ is a connected compact Hausdorff space, sometimes referred to as a \emph{base space}, and $r_n^i\in\mathbb{N}$, $p_n^i\in C(Z_n^i)\otimes M_{r_n^i}$ is a projection, for $i\in \{1,2,\ldots,N_n\}$. The C$^*$-algebras $A_n$ will be referred to as the \emph{building blocks} of the inductive limit. The direct summands of such building blocks are precisely the homogeneous C$^*$-algebras, and we will denote by \textbf{H} the class consisting of all finite direct sums of homogeneous C$^*$-algebras (i.e. C$^*$-algebras that have the same form as $A_n$). The $*$-homomorphisms $\phi_n$, which will be referred to as the \emph{connecting maps} of the inductive limit, are assumed to be unital and injective (note that by \cite[Theorem~2.1]{EGL05} we may always reduce to the injective case), and to be of a \emph{generalized diagonal} form: letting $\phi_n^{ji}$ denote the canonical composition 
\begin{equation*}
    A_n^i \hookrightarrow A_n \xrightarrow{\phi_n} A_{n+1} \twoheadrightarrow A_{n+1}^j,
\end{equation*}
we require that, for $a\in A_{n}^i$,
\begin{equation}\label{equation:phiIJN}
    \phi_n^{ji}(a)=\sum\limits_{y \in \mathcal{Y}(n,(j,i))}(a\circ\lambda_y)\otimes q_y
\end{equation}
where $\mathcal{Y}(n,(j,i))$ is some indexing set, $\lambda_y:Z_{n+1}^j\rightarrow Z_n^i$ is a continuous map and $q_y$ is either 0 or a projection in $C(Z_{n+1}^j)\otimes M_{s_n^{ji}}$ that is a line bundle (i.e. a complex vector bundle of dimension 1). Here we are imposing $r_{n+1}^j=\sum\limits_{i=1}^{N_n}r_n^is_n^{ji}$.

We also assume that $p_n^i=\sum\limits_{r=1}^{\mathrm{rank}(p_n^i)}\oplus p_{n,r}^i$ where each $p^i_{n,r}$ is a line bundle. The reason behind this assumption is to be able to use \cite[Lemma~6.1]{Li18} to obtain Cartan subalgebras of direct sums of homogeneous C$^*$-algebras. Note that our setting is equivalent to the setting in \cite[Section~4]{LR22}, except that here the notation emphasizes the way in which the summands of the building blocks embed into each other.

As in \cite[Section~4]{LR22}, the following are examples of C$^*$-algebras that fit naturally into our setting:

\begin{itemize}
\item AI-algebras (see \cite[Section~3]{klaus2}).
\item Villadsen algebras of the first kind \cite{Vil98}.
\item Tom's examples of non-classifiable C$^*$-algebras \cite{Toms}.
\item Goodearl algebras \cite{Goo} (see also \cite[Example~3.1.7]{Ror}).
\item AH-algebras models for dynamical systems in \cite[Example~2.5]{Niu}
\item Villadsen algebras of the second kind \cite{Vil99}.
\end{itemize}

Notice that in some of the above examples, the connecting maps are unitary conjugates of generalized diagonal maps, however we might as well assume there are no unitaries present as we still obtain the same inductive limit up to isomorphism (see for instance \cite[2.1]{ellrr0}).

\subsection{The Associated C$^*$-diagonal}\label{sec:ass C* diag}

In \cite[Section~4]{LR22} a C$^*$-diagonal was obtained for AH-algebras with generalized diagonal connecting maps. We briefly recall this construction, which is based on the construction in Section 6 of \cite{Li18}. Whilst the full construction of a groupoid and twist is not needed for our paper, we will still present it here for the sake of completeness.

Since the $p_{n,r}^i$ are line bundles, it follows by \cite[Lemma~6.1]{Li18} that\newline $C_n^i=\bigoplus\limits_{r=1}^{\mathrm{rank}(p_n^i)}p_{n,r}^i(C(Z_{n}^i)\otimes M_{r_n^i})p_{n,r}^i$ is a Cartan subalgebra of $A_n^i=p_{n}^i(C(Z_{n}^i)\otimes M_{r_n^i})p_{n}^i$. Let us now proceed to explain how to obtain the twisted groupoid model $(\mathcal{G}_n^i,\Sigma_n^ i)$ which is associated to the Cartan pair $(A_n^i,C_n^i)$. 

Let $\mathcal{G}_n^i=Z_{n}^i\times\mathcal{R}_{n}^i$ be the groupoid which is taken as a product of two groupoids $Z_{n}^i$ and ${\mathcal{R}}_{n}^i$ where the latter is the full equivalence relation on a set with $\mathrm{rank}(p_n^i)$ many elements. Since the $p_{n,r}^i$ are locally trivial line bundles, we can always find an open cover $\{V^{n,i}_{r,a}\}_a$ of $Z_n^i$ such that we can find partial isometries $v^{n,i}_{r,a}$ on $V^{n,i}_{r,a}$ whose domain projections are $e_{11}$ (canonical rank one projection) and whose range projections are $p_{n,r}^i(z)$, so that $v^{n,i}_{r,a}(z)=p^n_{i,r}(z)v^{n,i}_{r,a}(z)e_{11}$. From this we form the quotient space $\Sigma^{n,i}_{r,s}:=\bigsqcup\limits_{c,a}(\mathbb{T} \times (V^{n,i}_{s,c} \cap V^{n,i}_{r,a}))/ \sim$ where the relation is given by $(t,v)\sim (t^\prime,v^\prime)$ if and only if $v=v^\prime$, and if $v\in V^{n,i}_{s,c} \cap V^{n,i}_{r,a}$, $v^\prime\in V^{n,i}_{s,d} \cap V^{n,i}_{r,b}$ then $t^\prime=v^{n,i}_{r,b}(v^{n,i}_{s,d})^*v^{n,i}_{s,c}(v^{n,i}_{r,a})^*t.$ The twist is then given by $\Sigma^i_n = \bigsqcup\limits_{r,s}\Sigma^{n,i}_{r,s}$ Notice that $\Sigma^{n,i}_{r,r}$ is just the trivial twist over $Z^i_n$, but in general, if the line bundles $p^i_{n,r}$ are not globally trivializable, we must connect trivial twists in this way using the local triviality property of the line bundles.

Multiplication in $\Sigma_n^i$ is defined as follows, $(([t,v],(r,s)),([t^\prime,v^\prime],(r^\prime,s^\prime)))\in {\Sigma_n^i}^{(2)}$ if and only if $v=v^\prime$ and $s=r^\prime$, and then the product is given by $([tt^\prime,v],(r,s^\prime))$. The surjection $\Sigma^i_n \twoheadrightarrow \mathcal{G}^i_n$ is given by $(([t,v],(r,s)) \to (v,(r,s))$.

Following the construction in \cite[Section~4]{LR22} we then obtain groupoid models for the homomorphisms $\phi_n^{ji}$. Indeed, recall that the indexing set for the sum appearing in \eqref{equation:phiIJN} is denoted by $\mathcal{Y}(n,(j,i).$ For $y\in \mathcal{Y}(n,(j,i))$ let $q_y$ be the associated line bundle appearing in \eqref{equation:phiIJN}, and let $p_{y}: (\mathcal{G}_{y},\Sigma_{y})\rightarrow (\mathcal{G}_{n}^i,\Sigma_{n}^i)$ denote the groupoid model for $a\to a\circ \lambda_y$ (an explicit description of $p_y$ can be obtained from Lemma 3.2 in \cite{BL17}). Then letting $$(H_n,T_n)=\bigsqcup\limits_{i,j,n,y=y(i,j,n)}(G_y,\Sigma_y)$$ and $(\mathcal{G}_n,\Sigma_n) \cong \bigsqcup\limits_{i,n}(\mathcal{G}_n^i,\Sigma_n^i)$ be the groupoid model for $A_n$ we obtain groupoid models for $A_n\rightarrow A_{n+1}$ via \begin{equation}\label{equation:disjointpy}
    (\mathcal{G}_n,\Sigma_n) \xleftarrow{\bigsqcup p_y} (H_n,T_n) \hookrightarrow (\mathcal{G}_{n+1},\Sigma_{n+1})
\end{equation}where the first arrow is a surjective, proper and fibrewise bijective groupoid homomorphism and the second arrow is an injective groupoid homomorphism with open image. These give rise to the connecting map $\phi_n$ and then by \cite[Theorem~3.6]{BL17} we obtain a C$^*$-diagonal for our AH-algebra $A$, which we will label $C_A$ and call the \emph{associated C$^*$-diagonal} for $A$. Furthermore, \cite[Theorem~3.6]{BL17} describes the groupoid model for $(A,C_A)$.

\section{The Associated Bratteli Diagram and the Spectrum of the Associated C$^*$-diagonal}

For a pair $(A,C_A)$ as constructed above we will associate a certain Bratteli diagram $\mathcal{B}_A$ that will serve as a visual tool that will help us understand the spectrum of $C_A$. For a summand $A_n^i$ of $A_n$, once again consider its embedding into $A_{n+1}^j$ via $\phi^{ji}_n$ as in \eqref{equation:phiIJN}. Let there be $N_n$ nodes at the $n^{\text{th}}$ stage of $\mathcal{B}_A$, one for each summand of $A_n$, and $N_{n+1}$ nodes at the $(n+1)^{\text{th}}$ stage of $\mathcal{B}_A$, one for each summand of $A_{n+1}$. These nodes should be associated to the underlying spaces giving rise to each summand, i.e. the $i^{\text{th}}$ node of the $n^{\text{th}}$ stage of $\mathcal{B}_A$ should be associated to $Z_n^i$. For each $y\in \mathcal{Y}(n,(j,i))$ draw $\mathrm{rank}(p_n^i)$ many arrows from the node representing $A_n^i$ to the node representing $A_{n+1}^j$, each arrow labeled by $\lambda_y$. In this way the arrows represent the eigenvalue functions $Z_{n+1}^j \to Z_n^i$. Doing this for all $n$, $i$ and $j$ will yield the Bratteli diagram we call $\mathcal{B}_A$. 

\bex\cite[Example~3.1.7~\text{(Goodearl Algebra)}]{Ror}\label{ex:goodearl}. Let $Z$ be a compact Hausdorff space and let $\{k_n\}_n$ be a sequence of natural numbers such that $k_n$ divides $k_{n+1}.$ Then let $r_n < \frac{k_{n+1}}{k_n}$ be a natural number, or 0, for each $n\in \mathbb{N}$. Let $z_{n,1},z_{n,2},\ldots,z_{n,r_n}$ be points in $Z$ (if $r_n>0$). Then the Goodearl algebra $A=\varinjlim(A_n,\phi_n)$ has building blocks $A_n=C(Z)\otimes M_{k_n}$ and the connecting maps $\phi_n$ are defined by $$\phi_n(f)(z)=\mathrm{diag}(f(z_{n,1}),\ldots,f(z_{n,r_n}),f(z),f(z),\ldots,f(z))$$where there are $\frac{k_{n+1}}{k_n}-r_n$ instances of $f(z)$ appearing in the diagonal. Hence every stage of the Bratteli diagram $\mathcal{B}_A$ would have a single node, and we would draw $k_n$ many arrows for each eigenvalue map $\delta_{{z_{n,i}}}$ from the $n^{\text{th}}$ node into the $(n+1)^{\text{th}}$ node, as well as $k_n$ many arrows for each of the $\frac{k_{n+1}}{k_n}-r_n$ identity maps that appear as eigenvalue functions. 
\eex

\bex\cite[Example~2.5~\text{(AH-model for dynamical system)}]{Niu}\label{ex:min dyn}
Let $Z$ be a compact Hausdorff space, $\sigma$ a homeomorphism on $Z$ and consider the dynamical system $(Z,\sigma)$. The AH-model for it is given by $A=\varinjlim(A_n,\phi_n)$ where the building blocks are given by $A_n=C(Z)\otimes M_{2^{n-1}}$, and the connecting maps are given by $$\phi_n(f)=\mathrm{diag}(f,f\circ \sigma).$$Hence $\mathcal{B}_A$ contains one node at each stage and we draw, between the $n^{\text{th}}$ node and $(n+1)^{\text{th}}$ node, $2^{n-1}$ arrows for the eigenvalue function given by the identity map, and $2^{n-1}$ arrows for the eigenvalue function given by $\sigma$.
\eex

For $y\in\mathcal{Y}(n,(j,i))$ we can explicitly describe the map $p_y:(\mathcal{G}_{y},\Sigma_{y})\rightarrow (\mathcal{G}_{n}^i,\Sigma_{n}^i)$ as it is a groupoid model for $a\to a\circ \lambda_y$. In fact for our purposes we need to only do this for the unit spaces of the twisted groupoids.  

\blemma\label{lemma:unit space maps}
Let $ \mathcal{Y}(n,(j,i))$ be the indexing set associated to $\phi_n^{ji}$ in \eqref{equation:phiIJN}. Let $p_y:(\mathcal{G}_{y},\Sigma_{y})\rightarrow (\mathcal{G}_{n}^i,\Sigma_{n}^i)$ be the groupoid model associated to the map $a\to a\circ \lambda_y$ from \eqref{equation:phiIJN}, $y\in \mathcal{Y}(n,(j,i))$ (c.f. Lemma 3.2 in \cite{BL17} for an explicit description of $p_y$). Let $([1,x],(r,r))\in \Sigma_y^0$ (see the notation in Section \ref{sec:ass C* diag} above), where $x\in Z_{n+1}^j$ and $r$ corresponds to the summand of the image of $\phi_n^{ji}$ in the $(p^i_{n,r}\circ \lambda_y)^{\text{th}}$ corner. Then $p_y((([1,x]),(r,r)))=([1,\lambda_y(x)],(r,r))$.  
\elemma

\bproof
For an element $c \in C^i_n$, let $\tilde{c}$ denotes its identification in the corresponding reduced twisted groupoid C$^*$-algebra. Then as $p_y$ is the groupoid model for $c\to c\circ \lambda_y$ we have that $\widetilde{c\circ \lambda_y}(([1,x],(r,r)))=\widetilde{c}(p_y(([1,x],(r,r))))$. The result follows by \cite[Lemma~6.2]{Li18}.
\eproof 

\bdefin
For a Cartan pair $(A,C_A)$ as above and associated Bratteli diagram $\mathcal{B}_A$, a sequence $p=(\lambda_{y_1},\lambda_{y_2},\ldots)$, for $y_k\in \mathcal{Y}(k,(j,i))$ for some $i$ and $j$, will be called an \emph{infinite path} (or \emph{path} if the context is clear) if the source node of the edge with label $\lambda_{y_k}$ is the range node of the edge with label $\lambda_{y_{k-1}}$. The domain of a map $\lambda_y$ for $y\in \mathcal{Y}(n,(j,i))$ for some $n$, $i$ and $j$ will be denoted $D(y)$, and the codomain $C(y).$ For  a path $p$, we will denote by $\varprojlim p$  or $\varprojlim\limits_k\{C(y_k),\lambda_{y_k}\}$ the inverse limit of the sequence of spaces associated to the maps $\{\lambda_{y_k}\}$ with respect to the path $p$. The set of all paths will be denoted $P_{\mathcal{B}_A}$. A \emph{cylinder path} $Q\subseteq P_{\mathcal{B}_A}$ is a collection containing all paths which share the same initial $k$-edges for some $k\in \mathbb{N}$.
\edefin

\btheo\label{theorem:connected components}
For a Cartan pair $(A,C_A)$ as above and associated Bratteli diagram $\mathcal{B}_A$, the spectrum of $C_A$ is homeomorphic to $\bigsqcup\limits_{p\in P_{\mathcal{B}_A}}\varprojlim p$. The connected components of the spectrum of $C_A$ are exactly the inverse limits $\varprojlim p$ for a path $p$.
\etheo

\bproof
First note that the map $([1,x],(r,r))\to x$ for fixed $r$ and $x \in Z_{n+1}^j$ for some $n$ and $j$ defines a homeomorphism onto $Z_{n+1}^j$. Since we assume that our connecting maps $\phi_n$ are unital and injective generalized diagonal maps, it is clear that they map full elements to full elements, and hence by \cite[Remark~5.6]{Li18} the spectrum of $C_A$ is given by $\varprojlim\limits_n\{\mathcal{G}_n^0,\bigsqcup p_y\}$ where $\bigsqcup p_y$ maps $\mathcal{G}_{n+1}^0$ onto $\mathcal{G}_n^0$, as in \eqref{equation:disjointpy}. Hence by Lemma \ref{lemma:unit space maps} every element of  $\varprojlim\limits_n\{\mathcal{G}_n^0,\bigsqcup p_y\}$ can be identified with an element of $\varprojlim p$ for some path $p$ in $\mathcal{B}_A$. The union under all such paths $\bigsqcup\limits_p\varprojlim p$ then yields $\varprojlim\limits_n\{\mathcal{G}_n^0,\bigsqcup p_y\}=\mathrm{Spec}(C_A)$ where the topology is the one induced by the inverse limit topology. Under this topology, a basic open set is of the form $\bigsqcup\limits_{p\in Q\subseteq P_{\mathcal{B}_A}}\varprojlim p$. However this basic open set is a union of other basic open sets over cylinder paths with longer initial shared edges. Since for any path $p$, $\varprojlim p$ is a Hausdorff continuum (see \cite[Theorem~117]{IM12}), these must be the connected components. 
\eproof

\bdefin\label{definition:unique Cartan}
For a class of $C^*$-algebras $\mathfrak{c}$, we will say that an inductive limit Cartan subalgebra $C=\varinjlim(C_n,\phi_n)$ inside an inductive limit C$^*$-algebra $A=\varinjlim(A_n,\phi_n)$, where $C_n,A_n \in \mathfrak{c}$, is $\mathfrak{c}$-\emph{unique} if for another inductive limit $C^*$-algebra $B=\varinjlim(B_n,\psi_n)\cong A$ and inductive limit Cartan subalgebra $D=\varinjlim(D_n,\psi_n)$, with $D_n,B_n\in \mathfrak{c}$, we have that there exists a $*$-isomorphism of the pair $(A,C)$ onto $(B,D)$. 
\edefin

\bdefin
We will say that an AH-algebra $A$ arising from generalized diagonal connecting maps has \emph{incomplete spectrum} or is \emph{spectrally incomplete} if there exists a path $(\lambda_{y_1},\lambda_{y_2},\ldots)\in P_{\mathcal{B}_A}$ and a Hausdorff continuum $Z=\varprojlim\limits_{k}(C(y_k),\gamma_{y_k})$ where $\gamma_{y_k}:C(y_{k+1})\rightarrow C(y_k)$ are continuous functions, such that no connected component of $\mathrm{Spec}(C_A)$ is homeomorphic to $Z.$ Such a space $Z$ will be called an \emph{unseen spectral component}. An AH-algebra that is not spectrally incomplete will be called \emph{spectrally complete} or said to have a \emph{complete spectrum}.
\edefin

\bprop\label{theorem:non-uniqueness}
Let $A$ be an AH-algebra arising from unital and injective generalized diagonal connecting maps, with associated C$^*$-diagonal $C_A$ and Bratteli diagram $\mathcal{B}_A$. Suppose that $A$ is spectrally incomplete, and such that for some unseen spectral component $Z=\varprojlim\limits_{k}(C(y_k),\gamma_{y_k})$ we have that the AH-algebra $B$ obtained by replacing the path $\{\lambda_{y_k}\}$ with $\{\gamma_{y_k}\}$ is isomorphic to $A$. Then $A$ does not admit an \textbf{H}-unique inductive limit Cartan subalgebra. 
\eprop

\bproof
By Theorem \ref{theorem:connected components}, one of the connected components of $C_B$ is $Z$. Hence $\mathrm{Spec}(C_A)\ncong \mathrm{Spec}(C_B)$ and so $C_A\ncong C_B$. 
\eproof

\bremark
The reason for stating Proposition \ref{theorem:non-uniqueness} and not a more general statement just assuming that the spectra of the associated C$^*$-diagonals are not homeomorphic is because we will use its hypotheses to obtain uniqueness and non-uniqueness of inductive limit Cartan subalgebras in certain classes of AH-algebras arising from generalized diagonal connecting maps. 
\eremark

\section{Examples: Spectrally Incomplete AI-algebras}
We will now consider AI-algebras, which are inductive limit C$^*$-algebras $A=\varinjlim(A_n,\phi_n)$ where the building blocks take the form $A_n=\bigoplus\limits_{i=1}^{N_n}A_n^i$ and where the summands are \emph{interval algebras}: $A_{n}^i=C([0,1])\otimes M_{r_n^i}$. The class of interval algebras will be denoted by \textbf{I}. We assume here that the connecting maps are injective and unital. By \cite[Theorem~3.1]{klaus2} we have that AI-algebras arise from generalized diagonal connecting maps, where the eigenvalue functions $\{\lambda_y\}$ are just some continuous maps $[0,1]\rightarrow [0,1]$ and the projections $\{q_y\}$ are just the canonical minimal diagonal projections. 

\btheo\label{theorem:AI non uniqueness}
Let $A=\varinjlim(A_n,\phi_n)$, where each $A_n$ is of the form $C([0,1])\otimes M_{r_n}$, be a simple and unital AI-algebra (viewed as arising from generalized injective and unital diagonal connecting maps). Assume $A$ is spectrally incomplete. Then $A$ does not admit \textbf{I}-unique inductive limit Cartan subalgebras.
\etheo

\bproof
First we show that $\{r_n\}$ contains an increasing subsequence. Assume for a contradiction that the sequence stabilizes. We may assume then that $A=\varinjlim\limits_{n}(C[0,1] \otimes M_N,\phi_{n})$ for a fixed $N \in \mathbb{N}$. Let $\mu_n$ denote the injection $C[0,1] \otimes M_N \hookrightarrow A$. We can write $\phi_n(a) = a \circ \lambda_n$ for some surjective eigenvalue function $\lambda_n:[0,1] \twoheadrightarrow [0,1]$ (surjectivity follows by injectivity of the connecting map). 

Let $I_1 = \{f \in C[0,1] \otimes M_N : f(0)=0\}$, which is a proper ideal of $C[0,1] \otimes M_N$. There exists $t_1$ such that $\lambda_1(t_1)=0$. Let $I_2 = \{f \in C[0,1] \otimes M_N : f(t_1)=0\}$. Repeat this process choosing $t_2$ such that $\lambda_2(t_2)=t_1$ and  $I_3 = \{f \in C[0,1] \otimes M_N : f(t_2)=0\}$ and so on. All the $I_n$'s are proper ideals and $\phi_{n}(I_n) \subseteq I_{n+1}$. Define $I = \overline{\bigcup\limits_{n=1}^\infty\mu_n(I_n)}.$ Then $I$ is an ideal in $A$. It is non-zero because each $\mu_n$ is injective, and it is not all of $A$, because $\mu_1(h)$, where $h$ is the constant matrix with each entry value $\frac{1}{2}$, is not close to any element in $\mu_n(I_n)$ for all $n \in \mathbb{N}$. Hence $I$ is a proper ideal which is a contradiction to simplicity.

Hence without loss of generality we may assume $A=\varinjlim\limits_{n}(C[0,1] \otimes M_{r_n},\phi_{n})$ with $r_n > 2^n$ and $k_{n} := \frac{r_{n+1}}{r_n} > 2^{n}$. Since $A$ is spectrally incomplete there exists an unseen spectral component $Z= \varprojlim([0,1], \gamma_{n})$. Let the set of eigenvalue functions corresponding to $\phi_{n}$ be $\mathcal{F}_{n} = \{\lambda^{n}_{1}, \ldots, \lambda^{n}_{k_{n}}\}$. Replace any one of the functions by $\gamma_n$. Replace two more functions by the functions $g(t)=\frac{t}{2}$ and $h(t)=\frac{t+1}{2}$, and keep whatever functions have not been replaced the same. Call this new set of functions $\mathcal{G}_{n} = \{w^{n}_1,\ldots, w^{n}_{k_n}\}$. 

Let $B=\varinjlim\limits_{n}(C[0,1] \otimes M_{r_n},\psi_{n})$ where the $\psi_{n}$'s are generalized diagonal connecting maps with associated eigenvalue functions the elements of $\mathcal{G}_{n}$. The existence of the eigenvalue functions $g$ and $h$ ensure injectivity of the connecting maps (as the union of these functions' images is $[0,1]$), and Lemma 1.2 in \cite{klaus3.5} ensures (by using the functions $g$ and $h$) that we get simplicity of $B$. 

We have the following commutative diagram: 

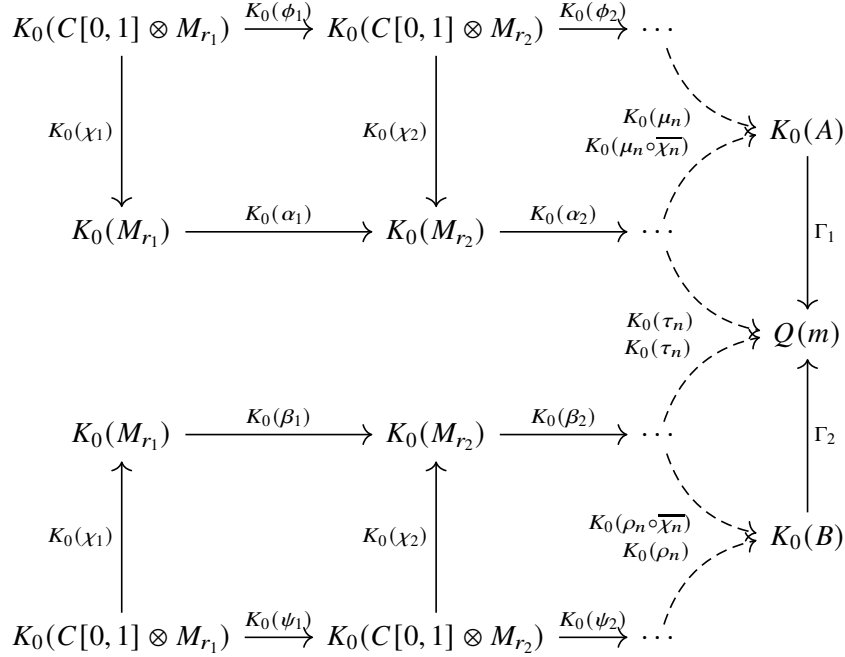
\begin{tikzcd}
{K_{0}(C[0,1] \otimes M_{r_1})} \arrow[r, "K_{0}(\phi_{1})"] \arrow[dd, "K_{0}(\chi_{1})"'] & {K_{0}(C[0,1] \otimes M_{r_2})} \arrow[r, "K_{0}(\phi_{2})"] \arrow[dd, "K_{0}(\chi_{2})"'] & \cdots \arrow[rd, "K_{0}(\mu_{n})"', dashed, bend right]                                                                                &                                  \\
                                                                                                &                                                                                                 &                                                                                                                                            & K_{0}(A) \arrow[dd, "\Gamma_1"]  \\
K_{0}(M_{r_1}) \arrow[r, "K_{0}(\alpha_{1})"]                                                 & K_{0}(M_{r_2}) \arrow[r, "K_{0}(\alpha_{2})"]                                                 & \cdots \arrow[ru, "K_{0}(\mu_{n} \circ \overline{\chi_{n}})", dashed, bend left] \arrow[rd, "K_{0}(\tau_{n})"', dashed, bend right] &                                  \\
                                                                                                &                                                                                                 &                                                                                                                                            & Q(m)                   \\
K_{0}(M_{r_1}) \arrow[r, "K_{0}(\beta_{1})"]                                                  & K_{0}(M_{r_2}) \arrow[r, "K_{0}(\beta_{2})"]                                                  & \cdots \arrow[ru, "K_{0}(\tau_{n})", dashed, bend left] \arrow[rd, "K_{0}(\rho_{n} \circ \overline{\chi_{n}})"', dashed, bend right]  &                                  \\
                                                                                                &                                                                                                 &                                                                                                                                            & K_{0}(B) \arrow[uu, "\Gamma_2"'] \\
{K_{0}(C[0,1] \otimes M_{r_1})} \arrow[r, "K_{0}(\psi_{1})"] \arrow[uu, "K_{0}(\chi_{1})"]  & {K_{0}(C[0,1] \otimes M_{r_2})} \arrow[r, "K_{0}(\psi_{2})"] \arrow[uu, "K_{0}(\chi_{2})"]  & \cdots \arrow[ru, "K_{0}(\rho_{n})", dashed, bend left]                                                                                   &                                 
\end{tikzcd}
\newline

In the diagram, $\mu_{n}: C[0,1] \otimes M_{r_n} \rightarrow A$ and $\rho_{n}: C[0,1] \otimes M_{r_n} \rightarrow B$  are the canonical injections of the building blocks into the respective inductive limits, $\chi_{n}$ is the map given by evaluation at 0, $\overline{\chi_{n}} : M_{r_n} \rightarrow C[0,1] \otimes M_{r_n}$ is the map sending $a$ to the continuous function with constant value $a$, and up to $K_0$, these are isomorphisms and inverses of each other. 

Indeed, for $s \in [0,1]$ the $*$-homomorphism $\phi_s: C[0,1] \otimes M_{r_n} \rightarrow C[0,1] \otimes M_{r_n} $ given by $\phi_s(f)(t)=f(st)$ defines a homotopy between $\overline{\chi_n} \circ \chi_n$ and $\mathrm{id}_{C[0,1] \otimes M_{r_n}}$. It is clear that $\chi_n \circ  \overline{\chi_n} = \mathrm{id}_{M_{r_n}}$.

The $\alpha_{n}$'s and $\beta_{{n}}$'s are the induced maps that make the diagram commutative, and $Q(m)$ is the subgroup of $\mathbb{Q}$ associated to the supernatural number $m$ corresponding to the inductive limit of the matrix algebras $M_{r_n}$'s, and the $\tau_{r_n}$'s are the normalized matrix traces. The isomorphisms $\Gamma_1$ and $\Gamma_2$ are canonically induced. Hence we get an isomorphism $$\phi_0 := \Gamma_2^{-1} \circ \Gamma_1 : K_{0}(A) \rightarrow K_{0}(B).$$ The diagram shows that $$\phi_{0}(K_{0}(\mu_{n})([p]_0 - [q]_0)) = [\rho_{n}(p)]_0 - [\rho_{n}(q)]_0,$$ and this is enough to determine the isomorphism. 

Now we determine an affine isomorphism $\phi_T: T_B \rightarrow T_A$, where $T_A$ ($T_B$) is the simplex of tracial states on $A$ ($B$). Note that the conditions of Lemma 4.1 in \cite{klaus3.5} are satisfied. Indeed, if $P_n$ is the center-valued trace appearing in that lemma, then $P_n \circ \phi_n - P_n \circ \psi_n$ has a factor $\frac{1}{r_{n+1}} < 2^{-(n+1)}$ appearing (as $P_n$ is normalized), and so the sum $\sum\limits_{n=1}^\infty \norm{P_n \circ \phi_n - P_n \circ \psi_n}$ is finite. Hence by Lemma 4.1 in \cite{klaus3.5} there is an affine isomorphism $\phi_T: T_B \rightarrow T_A$ satisfying, for all $i \in \mathbb{N}$, $\tau \in T_B$, and $a \in C[0,1] \otimes M_{n_i}$, $$\phi_T(\tau)(\mu_i(a)) = \mathrm{norm-lim}_{j \to \infty} \tau (\rho_j(\phi_{j,i}(a))).$$

Now note that for $g=K_{0}(\mu_{i})([p]_0 - [q]_0) \in K_0(A)$ and $\tau \in T_B$ we have that $$\langle \phi_0(g), \tau \rangle = \tau(\rho_{i}(p))-\tau(\rho_{i}(q))$$ (where the definition of $\langle \cdot,\cdot \rangle$ is given in Definition 1.1.10 in \cite{lars}). On the other hand, we have that \begin{equation}\label{longy}\langle g, \phi_T(\tau) \rangle = \mathrm{norm-lim}_{j \to \infty} \tau (\rho_j(\phi_{j,i}(p))) - \mathrm{norm-lim}_{j \to \infty} \tau (\rho_j(\phi_{j,i}(q))).\end{equation} Now note that we may choose the $p$ and $q$ in the definition of $g$ to be diagonals with entries either constant map 0 or constant map 1, as the trace of the projection determines the element in $K_0(A_i)^+$. Then $\phi_{j,i}$ and $\psi_{j,i}$ would agree on such projections and so we may replace $\phi_{j,i}$ in \eqref{longy} by $\psi_{j,i}$. Hence \eqref{longy} simplifies to $\tau(\rho_{i})(p) - \tau(\rho_{i})(q)$ and so we have $$\langle \phi_0(g), \tau \rangle = \langle g, \phi_T(\tau) \rangle.$$

By Theorem 2 in \cite{ell0}, we have that $$A \cong B.$$ By Theorem \ref{theorem:connected components} $Z$ is a connected component of $\mathrm{Spec}(C_B)$ but is an unseen spectral component of $\mathrm{Spec}(C_A)$. The result follows by Proposition \ref{theorem:non-uniqueness}.
\eproof

As an application of this result we may consider the class of Goodearl algebras that were considered in Example \ref{ex:goodearl}, restricted to unit interval base spaces. These form a particularly nice class of simple and unital AH-algebras arising from generalized diagonal connecting maps, some having real rank zero and some not. Consult \cite{Goo} or \cite[Example~3.1.7]{Ror} for further details. 

\bcor\label{cor: Goodearl}
Goodearl algebras with unit interval base spaces do not admit \textbf{I}-unique inductive limit Cartan subalgebras.
\ecor

\bproof
The eigenvalue functions that appear for Goodearl algebras are a combination of the identity function $f(x)=x$ and constant functions $f(x)=x_0\in[0,1]$. Hence given a path on the associated Bratteli diagram, it either contains infinitely many edges labeled with constant maps, or it does not. The first case would then correspond to a singleton inverse limit, whilst the second case would correspond to a Hausdorff continuum. Hence the connected components of the associated C$^*$-diagonal consists of Hausdorff continuua or singletons. However, there is a plethora of inverse limits on $[0,1]$ that are not of this form (see \cite[Chapter~1]{IM12}). Thus such Goodearl algebras are spectrally incomplete. The result follows by Theorem \ref{theorem:AI non uniqueness}.
\eproof

For another application we consider minimal dynamical systems $([0,1],\sigma)$ giving rise to an AH-algebra model as in Example \ref{ex:min dyn}. 

\bcor\label{cor: AH algebra models}
The AH-algebra model for a minimal dynamical system $([0,1],\sigma)$ does not admit \textbf{I}-unique inductive limit Cartan subalgebras.
\ecor

\bproof
The AH-algebra model is clearly unital arising from injective and unital generalized diagonal connecting maps, and because the dynamical system is minimal, it is simple. The eigenvalue functions on the edges of the corresponding Bratteli diagram are either the identity map $f(x)=x$ or the homeomorphism $\sigma$. Hence every inverse limit on a path is a Hausdorff continuum. Therefore the AH-algebra model is spectrally incomplete as there is a plethora of inverse limits on $[0,1]$ that are not Hausdorff continuua (see \cite[Chapter~1]{IM12}). The result follows by Theorem \ref{theorem:AI non uniqueness}.
\eproof

\section{Uniqueness of \textbf{AF}-Inductive Limit Cartan Subalgebras in AF-algebras}
In this section we consider the case of AF-algebras, which are spectrally complete. Perhaps unsurprisingly, they turn out to have \textbf{F}-unique inductive limit Cartan subalgebras, where \textbf{F} denotes the class of finite dimensional C$^*$-algebras. This result was already known, as a consequence of the work of Krieger in \cite{Kri}, and also \cite[Theorem~5.7]{Pow} provides a direct proof. Later in \cite{HP02}, the uniqueness of such masas is considered for more general limits. We wish to reprove this important result now, in a different way, as it has come to our attention that not many are aware of it. In fact our proof was written before we were made aware of these results. Such results would naturally be difficult to find for someone working with Cartan subalgebras as the terminology was different before the introduction of Cartan subalgebras for C$^*$-algebras. Our proof mimics Elliott's proof of the classification of AF-algebras via $K_0$ \cite{ellaf}, as presented in \cite[Chapter~7]{lars}, but now we keep track of Cartan subalgebras. Such techniques are more in line with the techniques used in Theorem \ref{theorem:AI non uniqueness}.

\blemma
AF-algebras are spectrally complete.
\elemma

\bproof
These algebras arise from generalized diagonal connecting maps with base spaces singletons. Any inverse limit over singleton sets is a singleton. 
\eproof

\blemma\label{Lemma 1}
Let $A=\bigoplus\limits_{j=1}^N M_{n_j}$ and $B=\bigoplus\limits_{i=1}^M M_{m_i}$ be finite dimensional $C^*$-algebras with Cartan subalgebras $C$ and $D$ respectively. Assume there exists an order unit preserving positive group homomorphism $\alpha: K_{0}(A) \rightarrow K_0(B)$. Then there exists a unital $*$-homomorphism $\phi : A \rightarrow B$ such that $K_0(\phi) = \alpha$, $\phi(C) \subseteq D$ and $\phi(N_{A}(C)) \subseteq N_{B}(D)$.
\elemma

\bproof
By Lemma 7.3.2 (i) in \cite{lars} there exists a unital $*$-homomorphism $\overline{\phi}: A \rightarrow B$ with $K_0(\overline{\phi})=\alpha$. Let $\{e_{pq}^j\}$ and $\{h_{uv}^i\}$ be systems of matrix units for $A$ and $B$ with respect to the Cartan subalgebras $C$ and $D$, respectively (obtained by extending the orthogonal basis for the Cartan subalgebras to a basis for the matrix algebras). Let $\{k_{ij}\}$ denote the multiplicity of the imbedding $M_{n_j}\hookrightarrow M_{m_i}$ via $\overline{\phi}$. Let us call this imbedding $\overline{\phi}_{ij}$, which is a $*$-homomorphism that is not necessarily unital. Since $\{\overline{\phi}_{ij}(e^j_{pp}) : 1 \le p \le n_j\}$ is a set of mutually orthogonal projections in $M_{m_i}$, each with trace $k_{ij}$, and $\sum\limits_{j=1}^N n_jk_{ij}=m_i$, we may decompose $\{h^i_{uu} : 1 \le u \le m_i\}$ into disjoint sets $H^{i1},H^{i2},\ldots,H^{iN}$ where each $H^{ij}$ has size $n_jk_{ij}$. Each $H^{ij}$ can then be decomposed into $n_j$ disjoint sets $H^{ij}_1,\ldots,H^{ij}_{n_j}$ of size $k_{ij}$ each. Let $f^{ij}_{pp} \in D$ be the sum of the elements of $H^{ij}_p$, with trace $k_{ij}$. Hence $\overline{\phi}_{ij}(e^j_{pp})$ has the same trace as $f^{ij}_{pp}$ and so $[\overline{\phi}_{ij}(e^j_{pp})]_0=[f^{ij}_{pp}]_0$.

Since $$f^{ij}_{11} \sim f^{ij}_{22} \sim \cdots \sim f^{ij}_{n_{j}},$$ we may apply Lemma 7.1.2 in \cite{lars} to obtain a system of matrix units $$\{f^{ij}_{pq} : p,q \in \{1,\ldots, n_j\}\}$$ in $M_{m_i}$. This can be extended to a set of matrix units $$\{f^{ij}_{pq} : j \in \{1,\ldots,N\}, p,q \in \{1,\ldots, n_j\}\}$$ in $M_{m_i}$. In fact we can be a bit more specific. If \begin{equation}\label{orderborder1}H^{ij}_{p}=\{h^i_{u_1u_1},\ldots, h^i_{u_{k_{ij}}u_{k_{ij}}}\}, \; \; u_1 \le u_2 \le \ldots \le u_{k_{ij}},\end{equation} the sum of whose elements is $f^{ij}_{pp}$, and \begin{equation}\label{orderborder2}H^{ij}_{q}=\{h^i_{v_1v_1},\ldots, h^i_{v_{k_{ij}}v_{k_{ij}}}\}, \; \; v_1 \le v_2 \le \ldots \le v_{k_{ij}},\end{equation} the sum of whose elements is $f^{ij}_{qq}$, we define \begin{equation}\label{sharoota}f^{ij}_{pq}= \sum\limits_{s=1}^{k_{ij}} h^i_{u_sv_s}.\end{equation} Because we have fixed an order in \eqref{orderborder1} and \eqref{orderborder2}, we indeed get a well-defined system of matrix units $\{f^{ij}_{pq}\}$, with $$f^{ij_1}_{p_1q_1}f^{ij_2}_{p_2q_2}= 0 \; \; \text{unless} \; \; j_1=j_2, q_1=p_2 \; \; \text{in which case we get} \; \; f^{ij_1}_{p_1q_2}.$$ Define $$\phi: A \rightarrow B \; \; \text{by} \; \; \phi(e^j_{pq})= \sum\limits_{i=1}^M f^{ij}_{pq},$$ and extend it linearly to $A$. It is clear that $\phi$ is a unital $*$-homomorphism with $\phi(C) \subseteq D$. The elements of $N_A(C)$ are those that have at most one non-zero entry in any row or column with respect to $\{e^j_{pq}\}$, and the elements of $N_B(D)$ are those that have at most one non-zero entry in any row or column with respect to $\{h^{i}_{uv}\}$. Hence it follows by using \eqref{sharoota} that $\phi(N_A(C)) \subseteq N_{B}(D)$. Now note that \begin{equation*}
\begin{split}
    K_0(\phi)([e^j_{11}]_0)& =\sum\limits_{i=1}^M[f^{ij}_{11}]_0=\sum\limits_{i=1}^M[\overline{\phi}_{ij}(e^j_{11})]_0=[\overline{\phi}(e^j_{11})]_0=K_0(\overline{\phi})([e^j_{11}]_0)\\ & =\alpha([e^j_{11}]_0)
    \end{split}
\end{equation*} and since the set $\{[e^j_{11}]_0 : 1 \le j \le N\}$ generates $K_0(A)$, it follows that $K_0(\phi)=\alpha$ as desired.
\eproof

\blemma\label{Lemma 2}
Let $A$ and $B$ be finite dimensional $C^*$-algebras with Cartan subalgebras $C$ and $D$ respectively, as in Lemma \ref{Lemma 1}. Assume $\phi,\psi:A \rightarrow B$ are unital $*$-homomorphisms which map $C$ into $D$, $N_A(C)$ into $N_A(D)$, and such that $K_0(\phi) = K_0(\psi)$. Then there exists $U \in \mathcal{U}(B) \cap N_B(D)$ such that $\psi = \mathrm{Ad}(U) \circ \phi$.
\elemma

\bproof
Let $\{e_{pq}^j\}$ and $\{h_{uv}^i\}$ be systems of matrix units for $A$ and $B$ with respect to the Cartan subalgebras $C$ and $D$, respectively. It is easy to see that the assumption $K_0(\phi)=K_0(\psi)$ implies that $K_0(\phi_{ij})=K_0(\psi_{ij})$ for $i \in \{1,\ldots,M\}$ and $j \in \{1,\ldots,N\}$. Hence the traces of $\phi_{ij}(e^j_{pp})$ and $\psi_{ij}(e^j_{pp})$ are the same for all $p \in \{1,\ldots,n_j\}$. Let $v^{ij}_1$ be a partial isometry witnessing the Murray-von Neumann equivalence between $\psi_{ij}(e^j_{11})$ and $\phi_{ij}(e^j_{11})$. Since by assumption these initial and range projections belong to $D$, they are each a sum of elements of the form $$\psi_i(e^j_{11})=\sum\limits_{s=1}^{k_{ij}}h^i_{u_su_s}, \; \; \phi_i(e^j_{11})=\sum\limits_{s=1}^{k_{ij}}h^i_{v_sv_s}.$$ Hence we can be particular with our choice of $v^{ij}_{1}$ by declaring it $$v^{ij}_{1}= \sum\limits_{s=1}^{k_{ij}}h^i_{v_su_s}.$$ It follows that $v^{ij}_{1}$ belongs to $N_{B_i}(D_i)$ (where $B_i$ and $D_i$ denote the restriction to the $i^\text{th}$ summand of $B$ and $D$ respectively). Then $v^{ij}_p$ is constructed as $\psi_i(e^j_{p1})v^{ij}_1\phi_i(e^j_{1p})$ and so also belongs to $N_{B_i}(D_i)$ by the assumption of the lemma. From this $U_i$ is then constructed as $\sum\limits_{j,p}v^{ij}_p$ and $U$ as $\bigoplus\limits_{i=1}^M U_i$.

To see that $U$ is a normalizer, let $h_{uu}^i \in D$. Note that because $\phi$ is unital and maps $C$ into $D$ there must exists some $j \in \{1,\ldots,N\}$ and $p \in \{1,\ldots,n_j\}$ such that $h^i_{uu}$ appears as a summand of $\phi(e^j_{pp})$. Hence $Uh_{uu}^iU^*=v^{ij}_ph^i_{uu}(v^{ij}_p)^*$ and so belongs to $D_i$. Hence $UDU^* \subseteq D$. A similar argument done by replacing $\phi$ with $\psi$ shows that $U^*DU \subseteq D$. Hence $U \in N_B(D)$. By construction $U\phi(e^j_{pq})U^*=\psi(e^j_{pq})$ and so $\psi = \mathrm{Ad}(U) \circ \phi$. 
\eproof

\btheo\label{main theorem}
Let $A$ and $B$ be unital AF-algebras with \textbf{F}-inductive limit Cartan subalgebras $C$ and $D$, respectively. Assume there exists an isomorphism
$$\alpha: (K_{0}(A),K_{0}(A)^+,[1_A]_0) \rightarrow (K_{0}(B),K_{0}(B)^+,[1_B]_0).$$

Then there exists a $*$-isomorphism $\phi:A \rightarrow B$ such that $K_0(\phi)=\alpha$ and $\phi(C) = D$. 
\etheo

\bproof
Let $B_0 = C_0 = D_0 =\mathbb{C}$, $\phi_0(\lambda)=\lambda 1_{A_1}$, $\psi_0(\lambda)=\lambda 1_{B_1}$, $\mu_0(\lambda)=\lambda 1_A$ and $\rho_0(\lambda)= \lambda 1_B$. We will assume $A$, $B$, $C$ and $D$ all arise as inductive limits as follows respectively:
\begin{equation*}
\begin{tikzcd}
B_0 \arrow[r,"\phi_0"] & A_1 \arrow[r,"\phi_1"] & A_2 \arrow[r,"\phi_2"] & \ldots \arrow[r,"\mu_n",bend left=60] & A,
\end{tikzcd}
\end{equation*}
\begin{equation*}
\begin{tikzcd}
B_0 \arrow[r,"\psi_0"] & B_1 \arrow[r,"\psi_1"] & B_2 \arrow[r,"\psi_2"] & \ldots \arrow[r,"\rho_n",bend left=60] & B,
\end{tikzcd}
\end{equation*}
\begin{equation*}
\begin{tikzcd}
C_0 \arrow[r,"\phi_0"] & C_1 \arrow[r,"\phi_1"] & C_2 \arrow[r,"\phi_2"] & \ldots \arrow[r,"\mu_n",bend left=60] & C,
\end{tikzcd}
\end{equation*}
\begin{equation*}
\begin{tikzcd}
D_0 \arrow[r,"\psi_0"] & D_1 \arrow[r,"\psi_1"] & D_2 \arrow[r,"\psi_2"] & \ldots \arrow[r,"\rho_n",bend left=60] & D.
\end{tikzcd}
\end{equation*}
We may assume that all the $\phi_{n}$'s, $\psi_{n}$'s, $\mu_{n}$'s and $\rho_{n}$'s are injective and unital $*$-homomorphisms, for $n=1,2,\ldots$, and given by block-diagonal imbeddings with respect to matrix units chosen as an extension of an orthogonal basis for the Cartan subalgebras. Hence we may assume that $\phi_n(N_{A_n}(C_n)) \subseteq N_{A_{n+1}}(C_{n+1})$. The same holds for the connecting maps $\psi_n$.

We have the following commutative diagram with positive and order unit preserving group homomorphisms:

\begin{equation*}
\begin{tikzcd}
K_0(B_0) \arrow[rr,"K_0(\rho_0)"] \arrow[dr, "K_0(\phi_0)",swap] & &  K_0(B) \\
& K_0(A_1) \arrow[ur, "\alpha \circ K_0(\mu_1)",swap]  & \\
\end{tikzcd}
\end{equation*}

From Lemma 7.3.3 in \cite{lars} we obtain that there exists $m_1 \in \mathbb{N}$ and a positive order unit preserving group homomorphism $\alpha_1$ such that we have a commutative diagram:

\begin{equation*}
\begin{tikzcd}
K_0(B_0) \arrow[rr,"K_0(\rho_0)", bend left=70] \arrow[r,"K_0(\psi_{m_10})"] \arrow[dr, "K_0(\phi_0)",swap] & K_0(B_{m_1}) \arrow[r, "K_0(\rho_{m_1})"] &  K_0(B) \\
& K_0(A_1) \arrow[u,"\alpha_1"]\arrow[ur, "\alpha \circ K_0(\mu_1)",swap]  & \\
\end{tikzcd}
\end{equation*}

Next, consider the following commutative diagram with positive and order unit preserving homomorphisms:

\begin{equation*}
\begin{tikzcd}
K_0(A_1) \arrow[rr,"K_0(\mu_1)"] \arrow[dr, "\alpha_1", swap] & &  K_0(A) \\
& K_0(B_{m_1}) \arrow[ur, "\alpha^{-1} \circ K_0(\rho_{m_1})",swap]  & \\
\end{tikzcd}
\end{equation*}

There exists $n_1 \in \mathbb{N}$ and a positive order unit preserving group homomorphism $\beta_1$ making the following diagram commute:

\begin{equation*}
\begin{tikzcd}
K_0(A_1) \arrow[rr,"K_0(\mu_1)", bend left=70] \arrow[r,"K_0(\phi_{n_11})"] \arrow[dr, "\alpha_1",swap] & K_0(A_{n_1}) \arrow[r, "K_0(\mu_{n_1})"] &  K_0(A) \\
& K_0(B_{m_1}) \arrow[u,"\beta_1"]\arrow[ur, "\alpha^{-1} \circ K_0(\rho_{m_1})",swap]  & \\
\end{tikzcd}
\end{equation*}

Repeating, we obtain a commutative intertwining:

\begin{equation*}
\begin{tikzcd}
K_0(B_{m_1}) \arrow[rr,"K_0(\psi_{m_2m_1})"] \arrow[dr, "\beta_1",swap] & & K_{0}(B_{m_2}) \arrow[dr,"\beta_2"] & \ldots & \ldots & K_{0}(B) \\
 & K_0(A_{n_{1}}) \arrow[ur,"\alpha_2"] \arrow[rr,"K_0(\phi_{n_2n_1})"] & & K_0(A_{n_2}) & \ldots & K_0(A) \arrow[u,"\alpha",swap]
\end{tikzcd}    
\end{equation*}

where all homomorphisms are positive and order unit preserving. Since inductive limits do not change if one takes a subsequence of the original building blocks, we may relabel to assume

\begin{equation}\label{zigzag1}
\begin{tikzcd}
K_0(B_{1}) \arrow[rr,"K_0(\psi_1)"] \arrow[dr, "\beta_1",swap] & & K_{0}(B_{2}) \arrow[dr,"\beta_2"] & \ldots & \ldots & K_{0}(B) \\
 & K_0(A_{1}) \arrow[ur,"\alpha_1"] \arrow[rr,"K_0(\phi_1)"] & & K_0(A_{2}) & \ldots & K_0(A) \arrow[u,"\alpha",swap]
\end{tikzcd}    
\end{equation}

By applying Lemma \ref{Lemma 1} and Lemma \ref{Lemma 2}, we lift \eqref{zigzag1} to a commutative intertwining

\begin{equation}\label{zigzag2}
\begin{tikzcd}
B_{1} \arrow[rr,"\mathrm{Ad}(V_2)\psi_{1}"] \arrow[dr, "g_1",swap] & & B_{2} \arrow[dr,"g_2"] & \ldots & \ldots & \overline{B} \\
 & A_{1} \arrow[ur,"f_1"] \arrow[rr,"\mathrm{Ad}(U_2)\phi_{1}"] & & A_{2} & \ldots & \overline{A} \arrow[u,"\Phi",swap]
\end{tikzcd}    
\end{equation}

where all the diagonal $*$-homomorphisms in \eqref{zigzag2} are unital, map Cartan subalgebra into Cartan subalgebra, normalizer set into normalizer set, and induce the corresponding $K_0$ group homomorphisms in \eqref{zigzag1}, and where all the unitaries $\{U_n\}$ and $\{V_n\}$ for $n=2,3,\ldots$ satisfy $U_n \in N_{A_n}(C_n)$ and $V_n \in N_{B_n}(D_n)$. We also have an induced $*$-isomorphism $\Phi$ between the induced inductive limits (see for example, Exercise 6.8 in \cite{lars}).  

We extend \eqref{zigzag2} to a commutative diagram

\begin{equation}\label{zigzag3}
\begin{tikzcd}
B_1 \arrow[rr,"\psi_1"]& & B_2 \arrow[rr,"\psi_2"] & & B_3 & \ldots & \ldots & B \\
B_1 \arrow[u,"\mathrm{id}"] \arrow[rr,"\mathrm{Ad}(V_2)\psi_1"] \arrow[dr,"g_1"] & & B_2 \arrow[u,"\mathrm{Ad}(V_2^*)"] \arrow[rr,"\mathrm{Ad}(V_3)\psi_2"] \arrow[dr,"g_2"] & & B_3 \arrow[u,"\mathrm{Ad}(\psi_2(V_2^*))\mathrm{Ad}(V_3^*)"]  \arrow[dr,"g_3"]& \ldots & \ldots & \overline{B} \arrow[u,"F"] \\
& A_1 \arrow[ur,"f_1"] \arrow[rr,"\mathrm{Ad}(U_2)\phi_1"]& & A_2 \arrow[ur,"f_2"] \arrow[rr,"\mathrm{Ad}(U_3)\phi_2"] & & A_3 & \ldots & \overline{A} \arrow[u,"\Phi"]\\
& A_1 \arrow[u,"\mathrm{id}"] \arrow[rr,"\phi_1"] & & A_2 \arrow[rr,"\phi_2"] \arrow[u,"\mathrm{Ad}(U_2)"] & & A_3 \arrow[u,"\mathrm{Ad}(U_3)\mathrm{Ad}(\phi_2(U_2))"] & \ldots & A \arrow[u,"G"]
\end{tikzcd}    
\end{equation}
where we have induced $*$-isomorphisms $F$ and $G$. Define $$\phi=F \circ \Phi \circ G.$$ We need to check $K_0(\phi) = \alpha$ and $\phi(C) = D$. Note that $K_0(\phi)$ and $\alpha$ agree on $K_0(\mu_n)(K_0(A_n))$ by the commutativity of \eqref{zigzag1} and \eqref{zigzag3}. Since $$K_0(A) = \bigcup\limits_{n} K_0(\mu_n)(K_0(A_n))$$ we obtain that $K_0(\phi) = \alpha.$ Note that from the commutativity of \eqref{zigzag3} it follows that $\phi$ maps $\mu_n(C_n)$ into $\rho_{n+1}(D_{n+1}) \subseteq D$. Hence $\phi$ maps $C$ into $D$ but as it is a $*$-isomorphism, it follows that $\phi(C) \subseteq D$ is a masa in $B$, and so it follows that $\phi(C)=D$.
\eproof

\bcor\label{cor:AF uniqueness}
Unital AF-algebras have \textbf{F}-unique inductive limit Cartan subalgebras.
\ecor

\bproof
If $(A,C)$ and $(B,D)$ are inductive limit Cartan pairs with $A \cong B$, then the ordered $K_0$ group of $A$ is isomorphic to the ordered $K_0$ group of $B$, and so we may find by Theorem \ref{main theorem} an isomorphism $(A,C) \cong (B,D)$, which gives uniqueness by Definition \ref{definition:unique Cartan}.
\eproof

\bremark
It is easy to conclude from Theorem \ref{main theorem} that if $C$ and $D$ are \textbf{F}-inductive limit Cartan subalgebras of an AF-algebra $A$, then there is an approximately inner automorphism of $A$ that maps $C$ onto $D$. Indeed, consider the pairs $(A,C)$ and $(A,D)$. Use Theorem \ref{main theorem} with $B=A$ and $\alpha=\mathrm{id}.$ One obtains an isomorphism $\varphi: (A,C)\rightarrow(A,D)$ such that $K_0(\varphi)=\mathrm{id}.$ Now extend \cite[Lemma~7.3.2(ii)]{lars} to inductive limits of finite dimensional C$^*$-algebras to obtain that $\varphi$ is an approximately inner automorphism. Hence the term \emph{uniqueness} for inductive limit Cartan subalgebras in AF-algebras implies in particular uniqueness via approximate unitary equivalence. 
\eremark

\end{document}